# The Deeper Roles of Mathematics in Physical Laws

Kevin H. Knuth
Departments of Physics and Informatics
University at Albany (SUNY), Albany NY, USA

"Familiarity breeds the illusion of understanding"
- Anonymous

**Abstract**
Many have wondered how mathematics, which appears to be the result of both human creativity and human discovery, can possibly exhibit the degree of success and seemingly-universal applicability to quantifying the physical world as exemplified by the laws of physics. In this essay, I claim that much of the utility of mathematics arises from our choice of description of the physical world coupled with our desire to quantify it. This will be demonstrated in a practical sense by considering one of the most fundamental concepts of mathematics: additivity. This example will be used to show how many physical laws can be derived as constraint equations enforcing relevant symmetries in a sense that is far more fundamental than commonly appreciated.

## Introduction

Many have wondered how mathematics, which appears to be the result of both human creativity and human discovery, can possibly exhibit the degree of success and seemingly-universal applicability to quantifying the physical world as exemplified by the laws of physics (Wigner, 1960; Hamming, 1980). That is, if the laws of physics are taken as fundamental, then how can it be that mathematics, which is often developed as a creative act, can possibly be relevant to the physical world? Or is mathematics somehow more fundamental than physics, and are we instead discovering the laws of mathematics and identifying which ones apply to or perhaps underlie particular physical situations?

These questions reside alongside longstanding questions regarding the nature of physical law. For example, Isaac Newton introduced the idea that some physical laws can be derived from other more fundamental physical laws—often given the more distinct title of *principles*. Today, many of these principles take the form of conservation laws or symmetries, which helps to account for their universality. However, this raises many questions, such as whether there exists a unique minimal set of fundamental principles from which everything derives, or whether some physical laws are derivable and others are determined by chance or decree.

Despite this, we can already see that mathematics plays at least two roles. The first role is related to symmetries, which necessarily represent foundational concepts since they are not readily derivable from more fundamental concepts. The second role is related to calculation where equations are used to quantify physical phenomena. While some of the equations are known to be derivable from more fundamental principles, some have been adopted as foundational concepts in their own right, such as specific Lagrangians.  This provides an important clue, especially since we know from experience that many equations have been arrived at through educated guesses or in some cases even trial-and-error.

In this essay, I will argue that much of the mathematics that applies to physics arises from our choice of description of the physical world coupled with our desire to quantify it.  I will demonstrate this in a practical example by considering one of the most fundamental concepts of mathematics: additivity.  This explicit demonstration reveals that there are two distinct aspects to the role that mathematics plays.  The first aspect is related to ordering and associated symmetries, and the second aspect is related to quantification and the equations that enable one to quantify things.  The result is that to a great degree the mathematical equations that we consider to represent physical laws arise as constraint equations that enforce basic symmetries resulting from our chosen descriptions and desire to quantify.  While it may be that not all equations can be arrived at in this way, by showing that this is more widespread than commonly appreciated, I aim to suggest that this is the reason why mathematics is so effective.

## Questions

When I was a physics student, I was once troubled by a seemingly simple question.  I went around the department posing my question to fellow students as well as professors.  Most people seemed to feel that what I was asking was a silly question and some had quick answers.  However, despite these situations, I never received what I felt to be a satisfactory or insightful response.[1]  I had asked:

> "Why is it that when I take two pencils and add one pencil, I always get three pencils?  And when I take two pennies and add one penny, I always get three pennies, and so on with rocks and sticks and candy and monkeys and planets and stars.  Is this true by definition as in 2+1 defines 3?  Or is it an experimental fact so that at some point in the distant past this observation needed to be verified again and again?"

The fact that this simple question is focused on a very familiar problem obscures the fact that it embodies the very same wonderment that many of us possess when considering the "unreasonable effectiveness" of mathematics in physics.  I have always felt that if we can come to understand the answer to this simple question, then we ought to begin to be able to understand the role that mathematics plays in physics at a whole.[2]  On the other hand, if this question leads us to running in circles then we clearly have a long way to go.

Now some may scoff at this question, as many of my fellow students and professors did.  One common answer was that this is just measure theory: "You have sets where the cardinality of the set is a measure, and measures on sets sum when you combine them via set union." However, what many respondents failed to realize was that this is an *axiom* of measure theory, which means that it is *assumed*. As a physicist interested in the foundations, assumptions cause me concern. And despite its familiarity, this is a big assumption indeed! Furthermore, taking this to be the answer to my query simply transforms my question to: "Why must set measures be additive?"

---

[1] In preparing for this essay, I was pleased to find that Hamming had posed a similar question: "I have tried, with little success, to get some of my friends to understand my amazement that the abstraction of integers for counting is both possible and useful. Is it not remarkable that 6 sheep plus 7 sheep make 13 sheep; that 6 stones plus 7 stones make 13 stones? Is it not a miracle that the universe is so constructed that such a simple abstraction as a number is possible? To me this is one of the strongest examples of the unreasonable effectiveness of mathematics. Indeed, l find it both strange and unexplainable." (Hamming, 1980)

[2] Einstein expressed a similar sentiment about particle physics: "I would just like to know what an electron is."

The answer to this is not obvious[3], and it is central to the problem of the unreasonable effectiveness of mathematics. Clearly addition works, and this might suggest that the mathematics we use has been selected because it has worked every time we have employed it, which would make the use of addition an experimentally-observed result. Now, it could very well be that this is how this all came about historically. That is, it could be that the mathematics that humanity developed was selected (evolved) to work again and again without violation.

But why all the fuss; the original question is just about counting isn't it?
Could it really be the case that this goes deeper?

## Symmetries and Order

There are some important symmetries at play when I combine sets of pencils or pennies or monkeys or stars[4]. If I combine set $A$ with set $B$ to form the joint set $D$, the result will be the same as if I combine set $B$ with set $A$. This symmetry is called *commutativity*. Considering sets joined by set union, since $D = A \cup B$ and $D = B \cup A$, we can write

$$A \cup B = B \cup A. \qquad \textbf{Commutativity}$$

Now there is another important symmetry called *associativity* where the order in which one unites sets of objects also does not affect the final results. That is, in combining three sets $A$, $B$, and $C$, I could first combine $A$ with $B$ and then combine the result with $C$, or I could combine $A$ with the result of the combination of $B$ with $C$, and so on. In mathematical parlance, we have that

$$(A \cup B) \cup C = A \cup (B \cup C) = (A \cup C) \cup B, \qquad \textbf{Associativity}$$

so that we can just drop the parentheses and simply write $A \cup B \cup C$.

I feel relatively confident that these conceptual symmetries are, and need to be, experimentally observed when combining many things such as pencils and monkeys and so on. This is mainly because there are some examples where these symmetries are not experimentally observed. For example, if I combine a burning match with a paper napkin and then sometime later combine them with a glass of water, I get a very different result than if I were to first combine the napkin with the glass of water and

---

[3] This is especially difficult since additivity of set measure is assumed by mathematicians. That is where they decided to start. Where else would one start? This highlights one of the subtle and insidious difficulties faced by those in foundational studies. Once it is decided to adopt an assumption as a foundational construct it precludes the ability to delve deeper within that framework. Furthermore, it discourages others from doing so through what I call "The Curse of Familiarity" where one's familiarity with a problem fosters an illusion of understanding that blinds one from seeing subtle clues, hints, connections and/or difficulties. The parable of Newton and the apple is an example where Newton momentarily saw through the familiarity of falling apples to realize a connection with the falling of the Moon about the Earth.

[4] Finite sets are sufficient for our purposes here as we are not attempting to model an infinite world.

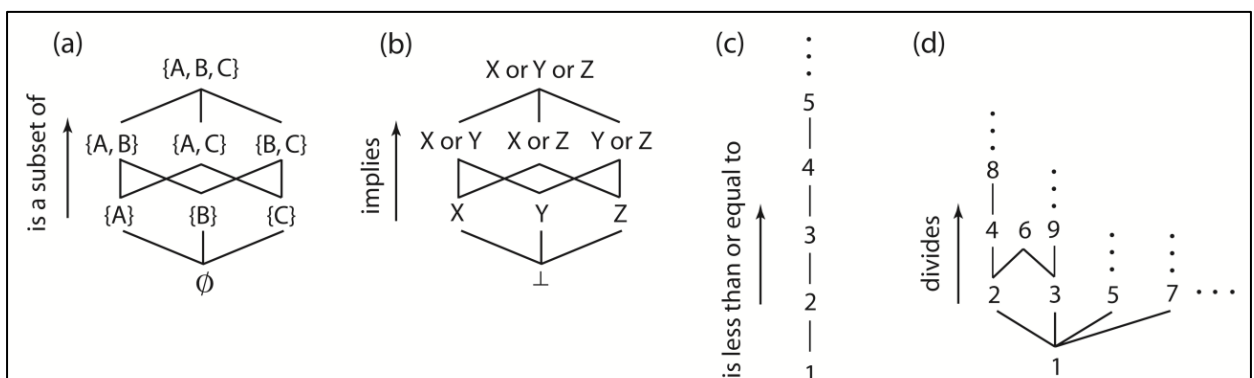

**Figure 1 (a)** The lattice of sets A, B, C ordered by set inclusion. The lattice meet and join are the operations set intersection and union, respectively. The bottom element is the null set. This structure is called a Boolean lattice as it represents a Boolean algebra. **(b)** The lattice of logical statements X, Y, Z ordered by implication also forms a Boolean lattice with the meet and the join being the logical 'and' and 'or' operations, respectively. The bottom element is the falsity. **(c)** The counting numbers ordered by the usual is-less-than-or-equal-to forms a lattice with the meet and join being the min and max functions. **(d)** However, the counting numbers can also be ordered according to whether one divides another with the meet and join being the GCD and LCM.

then combine those with the burning match. Of course, this depends upon our chosen level of description since if we included the air and counted the atoms, we would expect the results to be the same. This relies on the concept of closure where combining things results in the same sorts of things.

Sets are also interesting because they can be ordered. Sets can include other sets, and in fact we see this in the simple example above where the set $A \cup B \cup C$ includes the set $A \cup B$. In mathematical parlance, we say that the set $A \cup B$ is a subset of the set $A \cup B \cup C$ and write

$$A \cup B \subseteq A \cup B \cup C. \qquad \textbf{Subset Inclusion}$$

This leads to a hierarchy[5] of sets, which can be illustrated using what are called Hasse diagrams in order theory. Figure 1a illustrates this relationship among sets. The set union of two sets can be found by starting at each set and moving upward in the diagram until the two paths converge. Similarly, the set intersection can be found by starting at two sets and moving downward until the paths converge.

Order theory focuses on the ways sets of objects, called elements, can be ordered (Birkhoff, 1967; Davey & Priestley, 2002). The useful concept is that of a partially ordered set of elements. A *partially ordered set* (or *poset* for short) is defined as a set of elements that can be compared with a binary ordering relation, generically denoted $\leq$, that exhibits the properties of reflexivity, antisymmetry and transitivity (see Technical Endnotes).

The sets $A, B, C$ above along with the binary ordering relation of subset inclusion $\subseteq$ form a poset. However, this poset has additional properties that give it the special name of a lattice. Specifically, each pair of elements $x$ and $y$ has a unique least upper bound, $x \vee y$, called the join, and a unique greatest lower bound, $x \wedge y$, called the meet, such that the join and the meet are associative. So that in the case of our sets ordered by subset inclusion, we have that the lattice join is the set union, and the lattice meet is the set intersection.

---
[5] It would be better to call it a heterarchy since sets in general cannot be linearly ordered.

| (a) Sets, ⊆ | (c) Counting Numbers, ≤ |
|---|---|
| $A \subseteq B \iff \begin{array}{l} A \cap B = A \\ A \cup B = B \end{array}$ | $p \leq q \iff \begin{array}{l} \min(p, q) = p \\ \max(p, q) = q \end{array}$ |
| (b) Logical Statements, implies | (d) Counting Numbers, divides |
| $X \text{ implies } Y \iff \begin{array}{l} X \text{ and } Y = X \\ X \text{ or } Y = Y \end{array}$ | $p \text{ divides } q \iff \begin{array}{l} \text{GCD}(p, q) = p \\ \text{LCM}(p, q) = q \end{array}$ |

**Figure 2**. This figure illustrates how the ordering relation relates to the induced algebra for the four lattices illustrated in Figure 1. Later we will see that the fact that each of these algebras (set union and intersection, logical "and" and "or", min and max, and GCD and LCM) is associative, means that the same law (Sum Rule) will emerge when we consistently quantify the elements.

There are many other posets which form lattices. Figure 1 illustrates a few familiar examples. As described earlier, the first example in Figure 1a is the lattice of sets ordered by set inclusion. The second example (Figure 1b) shows the set of logical statements X, Y, Z ordered by logical implication where the lattice join is the logical *"or"* and the lattice meet is the logical *"and"*. Figure 1c, illustrates the counting numbers ordered by the usual less-than-or-equal-to ≤. This forms what is called a totally ordered set or a chain, where the lattice join and meet are the *max* and *min* functions, respectively. For example, $1 \vee 2 \equiv \max(1,2) = 2$ and $3 \wedge 4 \equiv \min(3,4) = 3$. Finally, Figure 1d illustrates that the counting numbers can be ordered in another way leading to a different poset. Here they are ordered according to whether one number divides another. The lattice join is given by the Least Common Multiple (*LCM*) and the lattice meet is the Greatest Common Divisor (*GCD*).

Now, what is important in all of this is that the simple act of ordering objects and the properties of commutativity and associativity allows one to view things in two different ways. On one hand we have a sort of hierarchy of elements determined by the binary ordering relation. And on the other hand, we can consider the join and meet to be algebraic functions that take two elements to a third. The result is that a lattice is an algebra, and the ordering relation is related to the algebraic relation by what is called the consistency relation

$$a \leq b \iff \begin{array}{l} a \wedge b = a \\ a \vee b = b \end{array} \quad \textbf{Consistency Relation}$$

With these concepts in mind, connections can now be made. Figure 2 shows the consistency relation written specifically for each of the four posets illustrated in Figure 1. One can now see that ordering, commutativity and associativity underlie a class of universal phenomena. I will next discuss how this leads to mathematics which gives rise to physical laws with a degree of universal applicability.

## Quantification

Let us now consider how we might go about quantifying elements of lattices. Two of the four examples above involve numbers directly and thus appear to be easy. There is nothing wrong with easy. Such cases will help us find the way since if there exists a general rule, it must apply to special cases.

We begin with the concept of quantification. The idea here is very simple. To each element $p$ we will assign a numeric value $v(p)$. Now, if we want our quantification scheme to maintain some representation of the ordering relation then for elements $p$ and $q$ where $p \geq q$ we assign values such that $v(p) \geq v(q)$. Essentially, here we are mapping elements to a total order—thus ranking them. This puts a strong constraint on the values we can assign.

Now given two disjoint elements $x$ and $y$ (which have null meet) and the values assigned to them $v(x)$ and $v(y)$, it would be good to know what number $v(x \vee y)$ we should assign to their join $x \vee y$. Clearly $v(x \vee y)$ must be greater than or equal to both $v(x)$ and $v(y)$. But can it be *any* such number? How can we be sure to avoid conflicts and ensure consistency with other joins in the same lattice?

First, if we want the assigned quantification to encode the underlying relationship, then we must assume that the number $v(x \vee y)$ we assign to the join $x \vee y$ of two disjoint elements, $x$ and $y$, is a function of the numbers $v(x)$ and $v(y)$. That is, we should be able to write

$$v(x \vee y) = v(x) \oplus v(y)$$

where $\oplus$ is a real-valued binary function to be determined. What could $\oplus$ possibly be? What would work and what would not work? For example, it couldn't be the min function since $\min(v(x), v(y))$ would result in picking the smallest number, which would violate the fact that the quantification must be greater than or equal to the largest number. The $\max$ function will not do either, since we will end up with a degenerate measure that simply ignores the lesser component, which is contrary to the goal of ranking.[6] The lesson is that not just any function will do.

Consider the symmetries of commutativity and associativity. First, since commutativity says that $x \vee y = y \vee x$ we have simply that their quantifications are equal $v(x \vee y) = v(y \vee x)$ which means that the operator $\oplus$ must be commutative: $v(x) \oplus v(y) = v(y) \oplus v(x)$. Writing $a = v(x)$ and $b = v(y)$, the commutative structure is more clear:

$$a \oplus b = b \oplus a. \qquad \textbf{Commutativity of } \oplus$$

Furthermore, since the lattice join is associative, the operator $\oplus$ must also be associative:

$$(a \oplus b) \oplus c = a \oplus (b \oplus c), \qquad \textbf{Associativity of } \oplus$$

where I have introduced a third element $z$ disjoint from $x$ and $y$, which has an associated value $c = v(z)$. The equation above is a functional equation called the Associativity Equation where the aim is to determine the function $\oplus$. The solution to this functional equation is known to be such that the function $\oplus$ must be an invertible transform of addition (Aczél, 1966; Craigen & Páles, 1989; Knuth & Skilling, 2012). That is $\oplus$ must take the form:

$$a \oplus b = f^{-1}(f(a) + f(b))$$

where the function $f$ is an arbitrary invertible function. In terms of the values we assigned to the elements this can be re-written as

$$f(v(x) \oplus v(y)) = f(v(x)) + f(v(x)).$$

---

[6] Since the goal is to rank elements via quantification (mapping elements to a total order), it may be helpful to formalize the desired systematic preservation of inequality by making explicit the assumption of cancellativity where for disjoint elements $x$, $y$, and $z$, where $v(x) \leq v(y) \leq v(z)$, we have $v(x) \oplus v(z) \leq v(y) \oplus v(z)$, which is implicit in any generally-useful (non-degenerate) notion of ranking.

This is significant because we can simply perform a *regraduation* on the valuations $v$ by instead assigning different values $u(x)$ defined by $u(x) = f(v(x))$ so that

$$u(x) \oplus u(y) = u(x) + u(y) \qquad \textbf{Additivity of } \oplus$$

is additive. The fact that we can *always* do this means that order, commutativity and associativity results in additive measures. We have, in fact, *derived* the countable additivity axiom of measure theory from a deeper symmetry principle!

## The Ubiquity of Additivity

We now can understand why we add quantities when we combine things. It doesn't always work. The properties of order, commutativity and associativity must apply to the selected description. So now when I take a pair of pencils and combine them with another pencil, I can quantify the union of these sets of pencils by simple addition. This works because sets of pencils are closed and can be ordered, and combining pencils is commutative and associative. This is the answer to my simple question in graduate school (and Hamming's question). And now we will see that the consequences are quite profound.

The additive rule above was derived for disjoint elements. One can show that in general (see Technical Endnotes) the additive rule is

$$u(x \vee y) = u(x) + u(y) - u(x \wedge y),$$

where we have to subtract off the value assigned to the meet of the two elements to avoid double-counting (Knuth, 2003; Knuth, 2010). In order theory, this is known as the *inclusion-exclusion principle* (Klain & Rota, 1997; Knuth, 2003). Others simply call it the *Sum Rule*. What is remarkable is that the Sum Rule appears over and over again, and now we can understand why. This is a consequence of closure, ordering, commutativity and associativity.

| Table 1  An Illustration of the Ubiquity of the Sum Rule | |
|---|---|
| $m(A \cup B) = m(A) + m(B) - m(A \cap B)$ | Measures on Sets |
| $P(A \text{ or } B \mid I) = P(A \mid I) + P(B \mid I) - P(A \text{ and } B \mid I)$ | Probability Theory |
| $\max(a, b) = a + b - \min(a, b)$ | Polya's Min-Max Rule |
| $\log(LCM(p, q)) = \log p + \log q - \log(GCD(p, q))$ | Integral Divisors |
| $I(A; B) = H(A) + H(B) - H(A, B)$ | Mutual Information |
| $\chi = V - E + F$ | Euler Characteristic |
| $E = (A + B + C) - \pi$ | Spherical Excess |
| $I_3(A, B, C) = \|A \sqcup B \sqcup C\| - \|A \sqcup B\| - \|A \sqcup C\| - \|B \sqcup C\| + \|A\| + \|B\| + \|C\|$ | Three Slit Problem |

Table 1 lists several examples of the Sum Rule in a variety of applications. The first four examples correspond to the four lattices in Figure 1. The example dealing with measures on sets is quite broad and includes physical volumes, surface areas, mean lengths, as well as linear superposition of potentials and other scalar quantities. The second example arises from the fact that probability quantifies the degree to which one logical statement implies another (Knuth & Skilling, 2012). This has enormous implications for physics since probability is one of the foundational concepts underlying both statistical mechanics and quantum mechanics. The third example, known as Polya's Min-Max Rule (Pölya & Szegö, 1964), is one of my favorites due to its simplicity and elegance. To find the larger of two numbers, sum them and take away the smallest! The fourth example illustrates that products are essentially linear (in the log) and they also arise from these basic symmetries (Knuth & Skilling, 2012).

Next, the definition of Mutual Information is a sum and difference of entropies, and is central to Information Theory (Cover & Thomas, 2012). This is followed by two less obvious examples related to geometry (Klain & Rota, 1997). The Euler Characteristic of a regular polytope is found by taking the number of Vertices minus the number of Edges plus the number of Faces. The spherical excess is an important quantity in spherical geometry where $A$, $B$, and $C$ are the angles of a spherical triangle.

The last example is the formula required to compute the quantum amplitude of the three slit problem where the square-cup operators represent combinations of slits (Sorkin, 1994). Here the quantity being summed is a complex number. Regardless of this, the formula still arises from associativity. It may be surprising, but one can *derive* the Feynman Rules for combining quantum amplitudes by relying on symmetries, such as associativity and distributivity along with consistency with Probability Theory (Goyal et al., 2010; Goyal & Knuth, 2011). And more recently, many of these same symmetries were used to derive the mathematics of flat spacetime as constraint equations describing the consistent quantification of a partially ordered set of events (Knuth, 2014; Knuth & Bahreyni, 2014; Knuth, 2015) as well as to derive the quantum symmetrization postulate (Goyal, 2015).

Now that we see how the quantification of phenomena is constrained by symmetries, and that this is far more widespread than the usual conservation laws and symmetries considered in physics, we can better understand precisely how mathematics is related to physical laws. The following questions now arise: "How far does this go?" and "To what degree are physical laws derivable and to what degree are they accidental, contingent or decreed by Mother Nature?"

## Conclusions

The results here shed light on the long-standing questions surrounding the unreasonable effectiveness of mathematics. We have seen interplay between two different aspects of mathematics. The first aspect is related to ordering and symmetries, and the second aspect is related to quantification and the equations that enable one to quantify things. Our choices in the phenomena that we focus on, the descriptions we adopt and the comparisons that we find important often amount to selecting a particular concept of ordering[7], which can possesses symmetries. The ordering relation and its symmetries in turn constrain consistent attempts at quantification resulting in constraint equations, which in many cases are related to what are considered to be physical laws. Much of the wonderment surrounding the unreasonable effectiveness of mathematics is not associated with the first aspect of ordering and symmetries since these more clearly depend on a choice of description and comparison, which in turn results in symmetries that can be easily observed and verified. Instead such wonderment is associated with the fact that we have equations that consistently allow us to quantify the physical world, and that these equations not only work very well, but in many cases exhibit some degree of universality. If we consider the equations themselves to be fundamental then the success of mathematics is somewhat of a mystery. But if we step back and release ourselves from familiarity and consider order and symmetry to be fundamental, then we see these equations as rules to constrain our artificial quantifications in accordance with the underlying order and symmetries of our chosen descriptions.

---

[7] Do we select a particular concept of ordering? Of course we do; an example is given in Figure 2 where the counting numbers are ordered in two different ways. Selecting one way of ordering gives you one set of laws (min and max) and selecting the other gives you another set of laws (GCD and LCM). The entire of field of number theory results as an attempt to study relationships between these two resulting sets of laws.

In some sense, this should not be surprising since it has been generally believed that the laws of physics reflect an underlying order in the universe. In fact, here it is explicitly demonstrated that some laws of physics not only reflect such order, but in fact can be derived directly from it. This has enormous implications for the direction and progress of foundational physics in the sense that it enables one to see that common mathematical assumptions, such as additivity, linearity, Hilbert spaces, etc., while familiar, are most likely not fundamental. Instead, there is room to delve deeper by identifying the fundamental symmetries and order-theoretic concepts that underlie physical theories.

The universality of mathematics, and specific mathematical relations, has been thought by many to be mysterious. Hamming wrote the following of Wigner's exposition (1960):

> "Wigner also observes that *the same mathematical concepts* turn up in entirely unexpected connections. For example, the trigonometric functions which occur in Ptolemy's astronomy turn out to be the functions which are invariant with respect to translation (time invariance). They are also the appropriate functions for linear systems. The enormous usefulness of the same pieces of mathematics in widely different situations has no rational explanation (as yet)." (Hamming, 1980)

As demonstrated in this essay, I believe that the answer lies in the deeper symmetries that various problems exhibit. To quote Jaynes:

> "the essential content ... does not lie in the equations; it lies in the ideas that lead to those equations." (Jaynes, 1959)

In addition to mystery and wonderment, this universality of mathematics has also caused a great deal of confusion. For example, some believe that the theory of quantum mechanics represents some kind of generalized or exotic probability theory. From the results presented above, one can see that while quantum amplitudes are used to compute probabilities, and while amplitudes and probabilities follow similar algebraic rules, quantum mechanics is not a generalized probability theory any more than information theory[8], geometry, and number theory are generalized probability theories. Instead, they all derive from the same fundamental symmetries. Each is a specific measure theory in its own right.

A deeper understanding of the roles that mathematics plays in physics in terms of order, symmetries, and quantification will help to clear up mysteries related to the effectiveness and universality of mathematics as well as to guide future foundational efforts.

---

[8] The entropies appearing in the definition of mutual information derive from probabilities, yet no one insists that probability theory is an exotic form of information theory. Entropy and probability are related in a very specific way with probabilities being used to compute entropies. Similarly, quantum amplitudes are related to probabilities in a very specific way with amplitudes being used to compute probabilities via the Born Rule.

# Technical Endnotes

This technical section provides the formal definitions of a partially-ordered set (poset) and a lattice as well as short derivation of the generalized sum rule for non-exclusive join.

A **partially-ordered set** (**poset**) $P$ is a set of elements $S$ along with a binary ordering relation, generically denoted $\leq$, postulated to have the following properties for elements $x, y, z \subseteq P$

| | |
|---|---|
| For all $x$, $x \leq x$ | Reflexivity |
| If $x \leq y$ and $y \leq x$, $x = y$ | Antisymmetry |
| If $x \leq y$ and $y \leq z$, $x \leq z$ | Transitivity |

A poset $P = (S, \leq)$ is referred to as a **partially** ordered set since not all elements are assumed to be comparable. That is, there may exist elements $x, y \subseteq P$ where it is neither true that $x \leq y$ nor that $y \leq x$. In these situations, we say that $x$ and $y$ are **incomparable**, which is denoted $x \parallel y$.

A **lattice** $L$ is a poset where each pair of elements has a supremum or least upper bound (LUB) called the **join**, and an infimum or greatest lower bound (GLB) called the **meet**. The meet of two elements $x, y \in L$ is denoted $x \wedge y$ and the join is denoted $x \vee y$. The meet and the join can be thought of as algebraic operators that take two lattice elements to a third lattice element. It is in this sense that every lattice is an algebra. The meet and join are assumed to obey the following relations

| | |
|---|---|
| $x \vee y = y \vee x$ <br> $x \wedge y = y \wedge x$ | Commutativity |
| $x \vee (y \vee z) = (x \vee y) \vee z$ <br> $x \wedge (y \wedge z) = (x \wedge y) \wedge z$ | Associativity |
| $x \vee (x \wedge y) = x$ <br> $x \wedge (x \vee y) = x$ | Absorption |

Lattice elements (and poset elements in general) can be quantified by assigning a real number (or more generally a set of real numbers) to each element. This is performed via a function $v$ called a valuation, which takes each lattice element to a real number: $v: x \in L \rightarrow \mathbb{R}$.

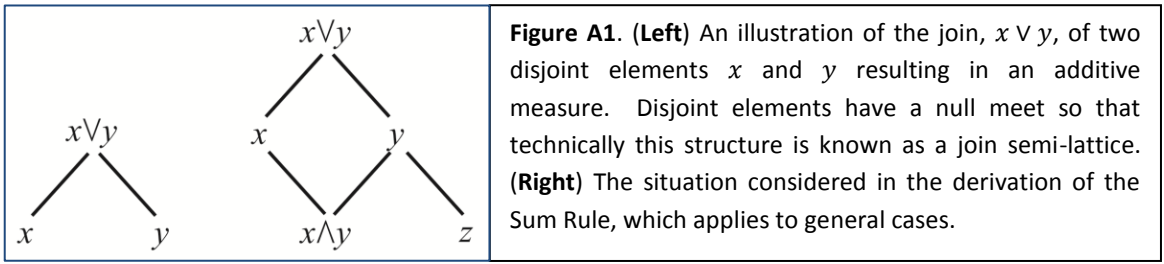

**Figure A1.** (**Left**) An illustration of the join, $x \vee y$, of two disjoint elements $x$ and $y$ resulting in an additive measure. Disjoint elements have a null meet so that technically this structure is known as a join semi-lattice. (**Right**) The situation considered in the derivation of the Sum Rule, which applies to general cases.

Valuations are meant to encode the ordering of elements in the lattice and this is accomplished by insisting that for $x \leq y$ we have that $v(x) \leq v(y)$. Furthermore, if the valuation is to encode the relationships among the elements, it must be that the valuation $v(x \vee y)$ assigned to the join of two disjoint elements $x$ and $y$ (Figure A1, Left) can be expressed as a function of the valuations $v(x)$ and $v(y)$ assigned to those two elements. We write this as

$$v(x \vee y) = v(x) \oplus v(y),$$

where the operator $\oplus$ is to be determined. The concept of a generally-useful quantification by valuation (non-degenerate ranking) implies that the operator $\oplus$ obeys a **cancellativity** property for disjoint elements $x$, $y$, and $z$, where for $v(x) \leq v(y) \leq v(z)$ we have $v(x) \oplus v(z) \leq v(y) \oplus v(z)$. Commutativity of the lattice join requires that the operator $\oplus$ is commutative, and associativity of the lattice join requires that the operator $\oplus$ is associative. The associative relationship represents a functional equation, known as the Associativity Equation, for the operator $\oplus$, whose solution is known to be an invertible transform of additivity (Aczél, 1966; Craigen & Páles, 1989; Knuth & Skilling, 2012), which can be written as

$$a \oplus b = f^{-1}(f(a) + f(b)),$$

where the function $f$ is an arbitrary invertible function. In terms of the valuations this is

$$f(v(x) \oplus v(y)) = f\big(v(x)\big) + f\big(v(x)\big).$$

This suggests that one can always choose a simpler quantification than the valuations $v$ by instead assigning values $u(x)$ defined by $u(x) = f(v(x))$ so that $\oplus$ is simple addition:

$$u(x) \oplus u(y) = u(x) + u(y).$$

Recall that this result holds only for disjoint elements (in a join semi-lattice). We now derive the result for two lattice elements in general. Consider the elements $x \wedge y$ and $z$ illustrated in Figure A1 (Right). Since their join is $y$, we have that

$$u(y) = u(x \wedge y) + u(z). \tag{A1}$$

Next consider that elements $x$ and $z$ are disjoint and their join is $x \vee y$. This allows us to write

$$u(x \vee y) = u(x) + u(z). \tag{A2}$$

Solving the (A1) for $u(z)$ and substituting into (A2) we have the Sum Rule

$$u(x \vee y) = u(x) + u(y) - u(x \wedge y), \qquad \text{Sum Rule}$$

which holds for general elements $x$ and $y$ (Knuth, 2003; Knuth, 2010; Knuth & Skilling, 2012).